\input amstex
\documentstyle{amsppt}
\input texdraw
\rightheadtext{Curvature, Diameter and Bounded Betti Numbers}
\leftheadtext{Zhongmin Shen and Jyh-Yang Wu }
\NoBlackBoxes
\magnification=1200


\define\e{\epsilon}

\topmatter
\title Curvature, Diameter and Bounded Betti Numbers
\endtitle
\author Zhongmin Shen and Jyh-Yang Wu
\endauthor
\date
Feb. 1999
\enddate
\address
\noindent
Z. Shen,
\newline
Department of Mathematical Sciences,
Indiana Univ.-Purdue Univ. at Indianapolis,
402 N. Blackford Street, Indianapolis, IN 46202-3216, U.S.A.
\newline
Email: zshen\@math.iupui.edu
\newline
\bigskip
J.-Y. Wu,
\newline
Department of Mathematics, National Chung Cheng University,
Ming-Hsiung, Chia-Yi 621, Taiwan
\newline
 Email: jywu\@math.ccu.edu.tw
\endaddress
\thanks J.-Y. Wu is partially supported by a Taiwan NSC grant.
\endthanks
\endtopmatter

\document

\head 1. Introduction\endhead

One of the fundamental problems in Riemannian geometry is to understand 
the relation between the topology and geometry of a Riemannian manifold. 

Every closed manifold $M$ admits  a Riemannian metric $g$
with the following curvature bound
$$
K  \geq - 1. \tag{1.1}
$$
Thus the curvature bound in (1.1) alone does not have any implication for the topological structure
of the manifold. With this normalized metric,
the topology depends on the ``size'' of the manifold. The diameter is one of the global geometric quantities to measure the manifold. Assume that 
$$
\text{ Diam} \leq D.\tag{1.2}
$$
It was proved by M. Gromov \cite{G1} that for any Riemannian $n$-manifold $(M, g)$ 
satisfying (1.1) and (1.2), 
the total Betti number (with respect to any field) is  bounded, namely, 
$$\sum_{i=0}^n \beta^i(M) \leq C(n, D).$$

\bigskip
In this paper, we are interested  in a much weaker curvature bound
$$
\text{ Ric} \geq -(n-1).\tag{1.3}
$$
First of all, according to Sha-Yang's examples \cite{SY},  the Gromov Betti Number Theorem 
is not true for Riemannian $n$-manifolds satisfying (1.2) and (1.3).
Nevertheless, the first Betti number is still bounded, i.e., 
$$ \beta^1 (M) \leq C(n, D).\tag{1.4}
$$
This is due to M. Gromov \cite{GLP} and T. Gallot \cite{GT}. 
 
Besides the ordinary  Betti number, what topological invariants are still controlled  by
the curvature bounds (1.1) or (1.3) ?
To answer this question, we
consider  the bounded 
cohomology groups $\hat{\text{H}}^*(M)$. The natural inclusion
$I: \hat{C}^*(M) \to C^*(M)$ induces a homomorphism
$\iota: \hat{\text{H}}^*(M) \to \text{H}^*(M)$ with image
$\tilde{\text{H}}^*(M):= \iota \Big [ \hat{\text{H}}^*(M) \Big ]$
(see Section 2 below for details).
 Put
$$
\align
 \tilde{\beta}^i(M) :& = \dim \tilde{\text{H}}^i (M),\\
\hat{\beta}^i(M) : & = \dim \hat{\text{H}}^i(M).
\endalign $$
Clearly,
$$ \tilde{\beta}^*(M) \leq \beta^*(M), \ \ \ \ \tilde{\beta}^*(M) \leq \hat{\beta}^*(M).$$
But there is no direct relationship  between $\beta^*(M)$ and $\hat{\beta}^*(M)$.
Both 
$\tilde{\beta}^i(M)$ and $\hat{\beta}^i(M)$ are  called the $i$th {\it bounded Betti numbers} of $M$. Below are some important examples:

\bigskip

(1)(Gromov)  For any closed manifold $M$,  $\hat{\text{H}}^1(M) =0$, hence
$\tilde{\beta}^1(M) = \hat{\beta}^1(M) =0$;

(2)(Thurston) For any closed manifold of negative curvature,  
$\tilde{\text{H}}^k(M) = \text{ H}^k(M) $ for all $k \geq 2$, hence
$\tilde{\beta}^k (M)= \beta^k(M)$ for $k \geq 2$ ;

(3) (Trauber) If $\pi_1(M)$ is amenable, then $\hat{\text{H}}^*(M) =0$, hence
$\tilde{\beta}^* (M) = \hat{\beta}^*(M)=0$ .

\bigskip

The bounded Betti numbers  behave quite differently from the ordinary Betti numbers. In particular, the Poincar\'{e} duality for $\beta^*(M)$ does not hold for
$\tilde{\beta}^*(M)$ and $\hat{\beta}^*(M)$.  By \cite{Gr3}, we know that the bounded Betti numbers $\hat{\beta}^*(M)$ are completely determined by $\pi_1(M)$. However, the bounded Betti numbers $\tilde{\beta}^*(M)$ depend not only on $\pi_1(M)$, but also on the higher dimensional topological structure of $M$, when $\pi_1(M)$ is not amenable. 

\bigskip

In this paper, we shall prove the following

\proclaim{Theorem A} There is a constant $C(n, D)$ depending only on $n, D$ such that for closed Riemannian $n$-manifold $(M,g)$ satisfying $\text{Ric}\geq -(n-1)$ and $\text{Diam}\leq D$,
 the total bounded Betti number
is bounded 
$$
\sum_{i=0}^n \tilde{\beta}^i (M) \leq C(n, D).\tag1.5
$$
\endproclaim

It is proved by Gromov that there is a small number $\e(n) >0$ such that if 
a closed oriented $n$-manifold
$M$  satisfies  (1.3) and $$
\sup_{p\in M}\text{vol} (B(p, 1)) < \e(n), \tag1.6
$$
 then there is a map 
$f$ from $M$ into   an $(n-1)$-dimensional polyhedron $P^{n-1}$ such that
the pre-image of any star neighborhood is contained in an {\it amenable} open subset. Then he concludes that $\|M\|=0$, which  is equivalent  to that $\tilde{\beta}^n (M)=0$.
Below is a generalization of Gromov's result.
\proclaim{Theorem B} Let $X$ be a compact metric space. There
is a small constant $\e = \e(n, X) >0$ such that for any 
 closed $n$-manifold $M$ satisfying $\text{Ric}\geq -(n-1)$, if 
$d_{GH}(M, X) < \e$, then 
$\tilde{\beta}^i(M) =0$ for $i > \dim X$. 
\endproclaim

\bigskip

Under a stronger curvature bound, the bounded 
cohomology groups $\hat{\text{H}}^*(M)$ of $M$  are controlled too. More precisely, 

\proclaim{Theorem C} There are only  finitely many isometric isomorphism
  types of bounded cohomology groups
 $(\hat{\text{H}}^*(M), \|\cdot\|_{\infty})$
among  closed Riemannian  $n$-manifolds satisfying 
$K\geq -1$ and $\text{Diam}\leq D$.
\endproclaim

From Theorem C, we conclude that there is a constant $C(n, D)$ such that 
 if a  closed $n$-manifold satisfying (1.1) and (1.2), then
$$ \hat{\beta}^*(M) \leq C(n, D), \tag1.7$$
provided that $\hat{\beta}^*(M) <\infty$. However, there are closed manifolds 
$M$ with $\hat{\beta}^*(M)=\infty$. Take
$M=S^1\times S^4 \#S^1\times S^4$. The fundamental group $\pi_1(M) = Z*Z$. Thus the second bounded cohomology
group $\hat{\text{H}}^2(M)$  is not finitely generated. This example is given to us by 
F. Fang.

\bigskip

\head 2. Preliminaries
\endhead
In this section, we shall summarize some of Gromov's results \cite{G3} which will be needed in our proof.

Let $M$ be a connected topological space. Denote by 
$\Sigma$ the set of all singular simplices $\sigma: \Delta \to M$.
The standard pseudo-$L^{\infty}$ norm on $C^*(M)$ is defined by
$$ \| c \|_{\infty}: =\sup_{\sigma \in \Sigma} |c(\sigma) |.$$
Consider the subcomplex
$\hat{C}^*(M)$ of bounded singular cochains of $M$. 
The homology groups $\hat{\text{H}}^*(M)$ of $\hat{C}^*(M)$ is called 
the {\it bounded cohomology} of $M$. 
Let $\|\cdot \|^b_{\infty}$ denote the induced norm on $\hat{\text{H}}^*(M)$.
Then $(\hat{\text{H}}^*(M), \|\cdot \|^b_{\infty})$ becomes a  normed space.
The natural inclusion
$I: \hat{C}^*(M) \to C^*(M)$ induces a homomorphism
$$ \iota: \hat{\text{H}}^*(M) \to \text{H}^*(M).\tag2.1$$
Put
$$ \tilde{\text{H}}^*(M): = \iota \Big [ \hat{\text{H}}^*(M)\Big ].$$
$\tilde{\text{H}}^*(M)$ is called  the {\it bounded part} of $\text{H}^*(M)$ (see \cite{Br}). 
Cohomology classes  in $\tilde{\text{H}}^*(M)$ are called {\it bounded classes} of $\text{H}^*(M)$.
Put
$$ \|\alpha \|_{\infty}: =\inf_{\beta\in I^{-1} (\alpha)}  \|\beta\|^b_{\infty}.\tag2.2$$
Then $(\tilde{\text{H}}^*(M), \|\cdot \|_{\infty})$ is a normed space.

\bigskip

  The most important fact in the bounded cohomology theory 
is that  the normed space
$(\hat{\text{H}}^*(M), \|\cdot \|_{\infty}^b)$ actually depends only on the fundamental group $\pi_1(M)$ (\cite{G3}) and  $\hat{\text{H}}^*(M) =0$ for connected closed manifolds with $\pi_1(M)$ amenable.

Below, we shall sketch Gromov's ideas to prove the above fact.
One is referred to \cite{Br}\cite{I} for different arguments.
 Although Gromov's theory is for general topological spaces, we shall focus on closed manifolds. 

\bigskip

First, Gromov introduced a notion of
multipcomplex.  A {\it (simplicial) multicomplex} is defined as a set $K$
divided
into the union of closed affine simplices $\Delta_{\sigma}\subset K$, $\sigma\in I$ such
that the intersection of any two simplices
$\Delta_{\sigma}\cap\Delta_{\tau}$ is a subcomplex in  both $\Delta_{\sigma}$ and $\Delta_{\tau}$.
The set $K$ with the weakest topology which agrees with the decomposition
$K=\cup_{\sigma\in I}\Delta_{\sigma}$ is denoted by $\vert K\vert$. The union of all
$i$-dimensional simplices in $K$ is called the {\it $i$-skeleton} of $K$,
denoted by $K^i\subset K$.

\bigskip

Let $M$ be an $n$-dimensional closed manifold. Denote by
$\Sigma$ the  set of all
singular simplicies $\sigma :\Delta^i\rightarrow M$, $i=0,1,2,\cdots$,
which are injective on the vertices of the standard (oriented)
$i$-simplex $\Delta ^i$. Take one copy of $\Delta^i$ for each $\sigma$, denoted by $\Delta^i_{\sigma}$, and put  $K:=\cup_{\sigma\in\Sigma}
\Delta^i_{\sigma}$. This union has a natural  structure of
a multicomplex such that
the canonical map:
$S: \vert K\vert\rightarrow M
$
defined by $S\vert_{\Delta^i_{\sigma}}:=
\sigma:\Delta^i_{\sigma}\rightarrow M$,  is
continuous.  Gromov proves that $S$ is a weak homotopy equivalence.
The multicomlex  $K$ is  
{\it large and complete} 
in the sense that  every component of $K$  has infinitely many
vertices, and  every
continuous map $f:\Delta^i\rightarrow K$ is homotopic, relative to $\partial \Delta^i$, to a simplicial
 embedding $g:\Delta^i\rightarrow K$,
provided $f\vert_{\partial \Delta^i}:\partial\Delta^i\rightarrow K$ is a simplicial
embedding.

\bigskip
For the  multicomplex $K$ constructed above,
there is   another   natural notion of  bounded cohomology $\hat{\text{H}}^*_a(K)$.
Let $\hat C^*_a(K)$ denote the complex of
bounded antisymmetric real cochains $c$, that is,
 $c(\Delta^i_{\sigma})=-c(\Delta^i_{\sigma\circ \delta})$
for any  orientation-preserving affine isomorphism 
$\delta: \Delta^i \to \Delta^i$. Then  $\hat{\text{H}}^*_a(K)$ is defined to be  the homology group
of $\hat C^*_a(K)$ with the natural 
 pseudo-norm $L_{\infty}$.
Define a homomorphism  $h: \hat{C}^i(M) \to \hat{C}^i_a(K)$ by
$$ h (\Delta^i_{\sigma}) := {1\over (i+1)!} \sum_{\delta} [\delta] c (\sigma \circ \delta ).\tag2.3 $$
Gromov asserts that $h$ induces an isometric isomorphism
$$ h^*: \hat{\text{H}}^i(M) \to \hat{\text{H}}^i_a(K).\tag2.4 $$

 In order to prove the fact that $\hat{\text{H}}^*(M)$ actually depend only on $\pi_1(M)$, 
Gromov introduces  a  large and complete subcomplex $i:\tilde{K} \hookrightarrow K$
with the following properties:

(i) each continuous map of a
simplex $\Delta^i$ into $\tilde{K}$,  whose restriction to the boundary is a
simplicial embedding,  is homotopic relative to the
boundary $\partial \Delta^i$ to
at most one simplicial embedding $\Delta^i\rightarrow K$.

(ii) the natural inclusion $i: \tilde{K} \hookrightarrow K$ is a homotopy equivalence. Hence it induces an isometric isomorphism
$$i^*: \hat{\text{H}}^*(K) \to \hat{\text{H}}^*_a (\tilde{K}). \tag2.5$$

\bigskip
A subcomplex $\tilde{K}$  with these properties exists and is  uniquely determined, up to an simplicial isomorphism, by the homotopy type of $K$. $\tilde{K}$ is called a {\it minimal model} of $K$.

\bigskip

Fix a minimal model $\tilde{K}$ of $K$. 
Let 
$\Gamma_1= \Gamma_1(\tilde{K})$ denote the group of simplicial automorphisms of $\tilde{K}$ which are homotopic to the identity
and keeps the $1$-skelton of $\tilde{K}$ fixed.
Then $\tilde{K}_1:=\tilde{K}/\Gamma_1$
 is a $K(\pi, 1)$ multicomplex 
with $\pi=\pi_1(\tilde{K})=\pi_1(M)$ and the 
projection $p: \tilde{K}\to \tilde{K}_1$ induces an isomorphism
between fundamental groups.  In particular,
the projection $p : \tilde{K}\rightarrow \tilde{K}_1$ induces an isometric isomorphism
$$ p^*: \hat{\text{H}}_a^*( \tilde{K}_1) \to \hat{\text{H}}_a^* ( \tilde{K}) .\tag2.6 $$

In virtue of (2.4)-(2.6), one can conclude that 
$$ \Phi: = p^{*-1}\circ i^* \circ h^*:
\hat{\text{H}}^*(M)\to \hat{\text{H}}^*_a(\tilde{K}_1)\tag2.7 $$
is an isometric isomorphism.
Thus Gromov concludes that 
  the normed cohomology groups  $\hat{\text{H}}^*(M)$   depend only on $\pi_1(M)$.

\bigskip

Let
$\tilde{\Gamma}:= \oplus_{x\in \tilde{K}^0_1}\pi_1(\tilde{K}_1, x)$.
The group  $\tilde{\Gamma}$ acts 
 on $\tilde{K}_1$ in a natural way.
Assume that  $\pi_1(M)$ is amenable, then $\tilde{\Gamma}$ is amenable.
The  standard averaging process 
leads to the following remarkable conclusion:
$\hat{\text{H}}^*(M) = \hat{\text{H}}^*_a(\tilde{K})=0$.
By a similar argument, one can show that   the amenable normal subgroups of $\pi_1(M)$ make no contributions to the bounded cohomology $\hat{\text{H}}^*(M)$. More precisely, we have the following

\proclaim{Lemma 2.1} Let $\Gamma \subset \pi_1(M)$ be a normal amenable subgroup. Then $\Gamma$ induces an action $\tilde{\Gamma}$ on $\tilde{K}_1$ such that $\tilde{K}_1/\tilde{\Gamma} $ is a multicomplex of $K(\pi, 1)$ type with $\pi = \pi_1(M)/\Gamma$ and 
$\hat{\text{H}}^*(M)$ is isometric isomorphic to
$\hat{\text{H}}^*_a(\tilde{K}_1/\tilde{\Gamma})$.
\endproclaim

Consider a class  ${\Cal M}$ of certain closed $n$-manifolds. 
Let
$$\align 
{\Cal M}_{\pi} : & =  \Big \{ \pi_1(M), \ M \in \Cal M\Big \}/\sim\\
\hat{\Cal M}^i : & = \Big \{ (\hat{\text{H}}^i(M), \|\cdot \|^b_{\infty}), \ M \in \Cal M\Big \}
\endalign $$
where $ \pi_1(M) \sim \pi_1(M')$ if and only if there are normal amenable subgroups $N \triangleleft \pi_1(M)$ and $N' \triangleleft \pi_1(M')$ such that
$ \pi_1(M)/N \approx \pi_1(M') /N'$.
Suppose that there are only finitely many isomorphism types of $\pi_1(M)$  in ${\Cal M}_{\pi}$. 
By Lemma 2.1, one can conlude that there are only finitely isometric isomorphism types of normed spaces $(\hat{\text{H}}^*(M), \cdot \|^b_{\infty})$
 in $\hat{\Cal M}^i$ for each $i$.

\bigskip

We now consider the bounded part $\tilde{\text{H}}^*(M)$ of $\text{H}^*(M)$. Althought 
  $\tilde{\text{H}}^*(M)$ is the image of 
 $\hat{\text{H}}^*(M)$, 
it is not clear 
how does  it depend on 
the fundamental group. In certain cases, the bounded cohomology group is very large, while 
the bounded part is trivial. Look at a closed integral homology $3$-spheres
$M$ with a hyperbolic metric. In this case,
$\tilde{\text{H}}^*(M)= \text{H}^*(M)=0$, but $\hat{\text{H}}^2(M) \not=0$.

It is natural to consider the case when a compact manifold $M$ is 
covered by a number of open amenable subsets. Here a subset $U$ is said to be amenable if for any $x\in U$, the image of the inclusion
$i_*: \pi_1(U, x) \to \pi_1(M, x)$ is an amenable subgroup.
One expects that $\tilde{\text{H}}^*(M)$ might be controlled by an amenable covering of the manifold.
Based on Gromov's bounded cohomology theory, 
N. V. Ivanov \cite{I} has made an important observation. He proved an analog of Leray's theorem on amenable coverings.

\proclaim{Lemma 2.2} \text{ (\cite{I})} Let $M$ be an $n$-dimrensional manifold , ${\Cal U}$ be an amenable covering of $M$, $N$ be the nerve of this covering, and $|N|$ be the geometric realization of the nerve. Then the canonical map
$\iota: \hat{\text{H}}^*(M) \to \text{H}^*(M)$ factors through the map
$\phi: \text{H}^*(|N|) \to \text{H}^*(M)$.  In other words,
there is a homomorphism $\psi:  \hat{\text{H}}^*(M)\to  \text{H}^*(|N|)$ such that $\iota = \phi \circ \psi$.
\endproclaim

\head 3. Proofs of Theorems A and  B \endhead 

Before we start to prove Theorem A, we recall a generalized
version of Margulis' lemma
which is due to Fukaya and Yamaguchi (\cite {FY1})
 in the sectional curvature case and then
extended by Cheeger and Colding (\cite {CC1}) to the Ricci curvature case.

\proclaim{Margulis' Lemma}
 Given $n$ and $k$, there exists a positive number $\epsilon (n)$
depending  only on $n$ and $k$
 such that if $M$ is a complete Riemannian $n$-manifold
with (1.3), then
 the image of the inclusion
 map $i_*:\pi_1(B(p,r),p)\rightarrow \pi_1(M)$ is almost
nilpotent for any point $p\in M$ when $r\le \epsilon(n)$.
\endproclaim

\noindent
{\it Proof of Theorem A}: 
Let $(M, g)$ be a closed Riemannian $n$-manifold satisfying
(1.2) and (1.3). Let $r = \e(n)$ be the number in 
Margulis's Lemma. Take a maximal set of disjoint
$r/2$-balls $B(p_i, r/2)$, $i=1, \cdots, m$. Then 
${\Cal U}:= \{ B(p_i, r)\}_{i=1}^m $ cover $M$. 
Assume that $B(p_{i_o}, r/2)$ has the smallest volume among $B(p_i, r/2)$.
By the Bishop-Gromov volume comparison, we obtain
$$ m \leq {\text{vol} (M) \over \text{vol}(B(p_{i_o}, r/2))} 
\leq {\int_0^D \sinh^{n-1}(t) dt \over \int_0^{r/2} \sinh^{n-1}(t) dt} = 
C(n,D).$$
Let $N$ be the nerve of this covering ${\Cal U}$  and $|N|$ be the geometric realization of the nerve. Since the number of the simplices in $N$ is bounded by $C(n, D)$, there is constant $C'(n, D)$ depending on $C(n, D)$ 
such that 
$$\dim \text{H}^*(|N|) \leq C'(n,D).$$
Note that each ball $B(x_i, r)$ in ${\Cal U}$ is amenable. 
By Lemma 2.2, we conclude that
$$ \tilde{\beta}^*(M) \leq \dim \text{H}^*( |N|) \leq C'(n, D).$$
This proves Theorem A.
\qed

\bigskip

The discussions above also suggest the following 
\proclaim{Problem}
Are there finitely many isometric isomorphism types 
of $(\tilde{\text{H}}^*(M), \|\cdot \|_{\infty} )$
among closed $n$-manifolds $(M, g)$
satisfying (1.2) and (1.3)?\endproclaim

Since there are infinitely many normed spaces in each dimension,
 Theorem A does not provide an answer to this question.

\bigskip

{\it Proof of Theorem B}. Let $X$ be as in Theorem B.
Since $X$ is compact, we can take
 a finite open covering $\{W_j\}$  of $X$ with mesh $<\epsilon (n)/8$
 and order $\le \dim X+1$,  where $\epsilon(n)$ is given by the Margulis'
Lemma.
 That is, $\text{Diam}(W_j)<\epsilon(n)/4$ for all $j$
 and every point $x$ is contained
no more than $\dim X+1$ subsets $W_j$.
\demo{Claim 1} There is a positive number $\delta$ such that every geodesic
ball $B(x,\delta)$ in $X$ is contained in some $W_j$.
\enddemo
Indeed, if this is not true, we can find a sequence of points $x_i$ in $X$ and
positive numbers
$\delta_i\rightarrow 0$
 such that the geodesic ball $B(x_i,\delta_i)$ is not totally
contained in any $W_j$ for all $i$. Since $X$ is compact, we can find, by
taking a subsequence if necessary, a limit point $x$ of $x_i$ in $X$. But now,
 the point $x$ must be in some $W_j$ and hence $W_j$ contains a geodesic ball
$B(x,r)$ for some positive raduis $r>0$. Then the triangle inequality implies
that   the geodesic ball $B(x_i,\delta_i)$ is contained in $W_j$ for large $i$.
This gives a contradiction and Claim 1 holds.

Next we consider the closed complement $F_j$ of $W_j$ in $X$,  
$F_j:=X-W_j$, and set
$$
E_j:=\Big \{x\in X\vert d(x,F_j)\ge \frac {\delta}2\Big \}.
$$
The set $E_j$ is closed and the triangle inequality implies that
$\{ E_j\}$ is a closed covering of $X$ due to our choice of $\delta$.

We take the positive number $\delta(X)$ to be the minimum of
$\delta/8$ and $\epsilon(n)/8$.

Assuming  the Gromov-Hausdorff distance between $M$ and $X$ is less than
$\delta(X)$, we can find an admissible metric $d$ on the disjoint
union $M\amalg X$ such that the classical Hausdorff distance of $M$ and $X$
in $M\amalg X$ is less than $\delta(X)$.
Then we define an open covering  ${\Cal U}= \{ U_j\}$ of $M$ by setting
$$
U_j:=\{ p\in M:d_{M\amalg X}(p,E_j)<2\delta (X)\}.
$$
The triangle inequality then gives that ${\Cal U}$ has mesh less than
$\epsilon(n)$ and it covers $M$.

\demo{Claim 2} The order of this open covering ${\Cal U}$
of $M$ is at most $\dim X+1$.
\enddemo
 Indeed, if there is a point $p$ in
$m= (\dim X+2)$ different open sets in ${\Cal U}$, say,
$U_j, j=1, \cdots, m$, 
 then we can find a point $x\in X$ with $d_{M\amalg X}(p,x)<\delta(X)$
 and the triangle inequality gives $d(x,E_j)<3\delta(X)$. Hence, one has
$d(p,F_j)\ge \frac {\delta}2-\frac {3\delta}8=\frac {\delta}8>0$ and thus
$p\in W_j$ for $j=1,2,\dots ,m$. This contradicts to the order of the covering
$\{ W_j\}$ since $m=\dim X+2$ and Claim 2 follows.

Therefore, we obtain an amenable open covering ${\Cal U}$ of $M$ when
$d_H(M,X)<\delta(X)$. Let $N$ be the nerve of this covering ${\Cal U}$ and $|N|$ be the geometric realization of the nerve. By our construction,
$$ \dim |N| \leq \dim X.$$
Thus
$$ \text{H}^i(|N|) = 0,  \ \ \ \ \ i > \dim X.$$
By Lemma 2.2, we conclude that  $\tilde{\beta}^i(M) =0$ for all
$ i > \dim X$.
\qed

 \head 4. Proof of Theorem C\endhead

In this section we shall prove  Theorem C.
First, we  recall the notion about the equivariant Hausdorff distance from
\cite {FY1}.
Let $\Cal M_{met}$ denote the set of all isometry classes of pointed inner metric
 spaces $(X,p)$ such that for each $r$ the ball $B(p,r)$ is relatively
 compact in $X$. Let $\Cal M_{eq}$ be the
 set of triples $(X,G,p)$ where $(X,p)$ is in $\Cal M_{met}$ and $G$ is a
 closed group of isometries of $X$. For $r >0$, put
$$
G(r)=\{g\in G\; \vert\;  d(gp,p)<r\}.
$$

\demo{Definition 4.1} Let $(X,G,x)$, $(Y,H,y)$ be in $\Cal M_{eq}$. An
 $\epsilon$-equivariant
pointed Hausdorff approximation stands for a triple $(f,\phi,\psi)$ of maps
$f:B(x,\frac 1{\epsilon})\rightarrow Y$, $\phi:G(\frac 1{\epsilon})\rightarrow
H(\frac 1{\epsilon})$ and $\psi:H(\frac 1{\epsilon})\rightarrow
G(\frac 1{\epsilon})$ such that
\roster
\item $f(x)=y$,
\item the $\epsilon$-neighborhood of $f(B(x,\frac 1{\epsilon}))$ contains
$B(y,\frac 1{\epsilon})$,
\item if $p, q\in B(x,\frac 1{\epsilon})$, then $\vert d(f(p),f(q))-d(p,q)\vert
<\epsilon$,
\item if $p\in B(x,\frac 1{\epsilon})$, $g\in G(\frac 1{\epsilon})$,
$gp\in B(x,\frac 1{\epsilon})$, then $d(f(gp),\phi(g)(f(p)))<\epsilon$,
\item if $p\in B(x,\frac 1{\epsilon})$, $h\in H(\frac 1{\epsilon})$,
$\psi(h)(p)\in B(x,\frac 1{\epsilon})$, then
$d(f(\psi(h)(p)),h(f(p)))<\epsilon$.
\endroster
\enddemo

We remark that it is required neither that $f$ is continuous nor that $\phi$,
$\psi$ are homomorphisms. The equivariant pointed Hausdorff distance
$d_{eH}((X,G,x),(Y,H,y))$ is defined to be  the infimum of the positive numbers
 $\epsilon$ such that there exist $\epsilon$ equivariant Hausdorff
 approximations from $(X,G,x)$ to $(Y,H,y)$ and from $(Y,H,y)$ to $(X,G,x)$.
By $d_{H}$ we denote the pointed Hausdorff distance, which is the case when
the groups are trivial. The notion
$$
\lim_{i\rightarrow\infty}(X_i,G_i,x_i)=(Y,G,y)
$$
means
$$
\lim_{i\rightarrow\infty}d_{eH}((X_i,G_i,x_i),(Y,H,y))=0.
$$

Now we proceed to prove Theorem C by the method of absurity as in \cite{W2}.
Suppose Theorem C were false. Then, there exists a sequence of Riemannian
$n$- manifolds $M_j$ satisfying $K\geq -1$ and $\text{Diam}\leq D$
 such that all of their bounded cohomology
$\hat{\text{H}}^*(M_j)$ are different.

Choose a base point $x_j$ in $M_j$ and a corresponding point $\tilde x_j$ in
its universal covering $\tilde M_j$. The fundamental group $\pi_1(M_j)$ acts
on $\tilde M_j$ as deck transformation.
Applying \cite {F} Theorem 2.1, \cite {FY1}
Proposition 3.6 and \cite {FY2} Theorem 4.1 for our sequence $(M_j,x_j)$
 and
their universal coverings $(\tilde M_j,\tilde x_j)$
 and fundamental groups $G_j=\pi_1(M_j)$,
one has
\proclaim {Lemma 4.1} There exist an Alexandrov space $(Y,y)$ and a Lie group
$G$ which is a closed subgroup of isometries of $Y$ such that one has
$Y/G=X$ and
$$
\lim_{i\rightarrow\infty}(\tilde M_i,G_i,\tilde x_i)=(Y,G,y).
$$
Moreover, for any  normal subgroup $H$ of $G$ with the properties
\roster
\item $G/H$ is discrete, and
\item $H$ is generated by $H(r)$ with $r>0$,
\endroster
 there exists a sequence of normal subgroups $H_i$ of $G_i$
such that
\roster
\item $\lim_{i\rightarrow\infty}(\tilde M_i,H_i,\tilde x_i)=(Y,H,y)$,
\item $G_i/H_i$ is isomorphic to $G/H$ for sufficiently large $i$,
\item $H_i$ is generated by $H_i(r+\epsilon_i)$ for some
$\epsilon_i$ with $\epsilon_i\rightarrow 0$.
\endroster
\endproclaim

Next we take the normal subgroup $G_0$ of the connected component of the
identity element of $G$. Since $G_0$ is generated by $G_0(\epsilon)$ for
any positive number $\epsilon$. We can choose $\epsilon$ to be
$\epsilon(n)/4$ where $\epsilon (n)$ is given by Margulis' Lemma.

Since $G$ is a Lie group, $G/G_0$ is discrete. Lemma 4.1 then implies that
there exists a sequence of normal subgroups $E_j$ of $G_j$ such that
$G_j/E_j$ is isomorphic to $G/G_0$ for sufficently large $i$. Moreover,
$E_j$ is generated by $E_j(2\epsilon)$ for large $j$.

From our choice of the number $\epsilon$, Margulis' Lemma implies that the
 normal subgroup $E_j$ is almost nilpotent.
Then, Lemma 2.1 yields that $\hat{\text{H}}^*(M_j)$ is isometrically isomorphic  to
$\hat{\text{H}}^*_a(K_j)$ for a multicomplex $K_j$ of $K(\pi, 1)$ type with
$\pi=G_j/E_j\simeq G/G_0$.

Since any two $K(\pi,1)$ multicomplices with isomorphic $\pi$'s are homotopy
equivalent, thus $\hat{\text{H}}^*(M_j)$ is isometrically
isomorphic
to each other for sufficiently large $j$. This leads to a contradiction and
Theorem C follows.
$\quad\square$

\enddocument